\documentclass[11pt]{article}
\usepackage{amsmath}
\usepackage{mathrsfs}
\usepackage{amsthm}
\usepackage{amsfonts}
\usepackage{amssymb}
\usepackage{bm}
\usepackage{mathtools}
\usepackage{color}
\usepackage{tikz}

\usepackage{enumerate}
\usepackage{enumitem}
\setlist[enumerate,1]{label=(\arabic*),font=\textup,
leftmargin=7mm,labelsep=1.5mm,topsep=0mm,itemsep=-0.8mm}
\setlist[enumerate,2]{label=(\alph*).,font=\textup,
leftmargin=7mm,labelsep=1.5mm,topsep=-0.8mm,itemsep=-0.8mm}

\usepackage[pdfstartview=FitH,bookmarksopen=false,
colorlinks=true,linkcolor=blue,citecolor=blue]{hyperref}

\newtheorem{theorem}{Theorem}[section]

\newtheorem{lemma}{Lemma}[section]

\title{\bf Spectral extremal graphs for edge blow-up of star forests \thanks {Research was partially supported by the National
Nature Science Foundation of China (grant numbers 11871329, 11971298, 1220010800) and Hainan Provincial Natural Science Foundation of China (No. 122QN218)}}
\author {Jing Wang$^{1}$  \, Zhenyu Ni$^{2}$ \, Liying Kang$^{3}$\thanks{\em Corresponding author. Email address: lykang@shu.edu.cn (L. Kang), wj517062214@163.com(J. Wang), 1051466287@qq.com(Z. Ni), fanyz@ahu.edu.cn(Y. Fan)} \,  Yi-zheng Fan\\
{\small $^{1}$College of Mathematics and Information Science, Henan Normal University,}
\\{\small Xinxiang 453007, P.R. China}\\
{\small $^{2}$Department of Mathematics, Hainan University,
Haikou 570228, P.R. China}\\
{\small $^{3}$Department of Mathematics, Shanghai University,
Shanghai 200444, P.R. China}\\
{\small $^{4}$ School of Mathematical Sciences, Anhui University,
Hefei 230601, P. R. China}}

\date{}
\begin{document}

\maketitle

\begin{abstract}
The edge blow-up of a graph $G$, denoted by $G^{p+1}$, is obtained by replacing each edge of $G$ with a clique of order $p+1$, where the new vertices of the cliques are all distinct. Yuan [J. Comb. Theory, Ser. B, 152 (2022) 379--398] determined the range of the Tur\'{a}n numbers for edge blow-up of all bipartite graphs and the exact Tur\'{a}n numbers for edge blow-up of all non-bipartite graphs. In this paper,
we   prove that the graphs with the
maximum spectral radius in an $n$-vertex graph without any copy of edge blow-up of star forests are the  extremal graphs for edge blow-up of star forests when $n$ is sufficiently large.


\medskip

\noindent{\bf AMS (2000) subject classification:}  05C50; 05C35
\end{abstract}

\section{Introduction}


Given a simple graph $G=(V(G),E(G))$, let $\delta(G)$ and $\Delta(G)$ be its minimum and maximum degrees of $G$, respectively. We use $e(G)$ to denote the number of edges in $G$. For $V_1,V_2 \subseteq V(G)$, $E(V_1,V_2)$ denotes the set of edges between $V_1$ and $V_2$, and $e(V_1,V_2)=|E(V_1,V_2)|$. For any $S\subseteq V(G)$, let $N(S)=\cup_{u\in S}\{v: uv\in E(G)\}$, $d_{S}(v)=|N_{S}(v)|=|N(v)\cap S|$. $G\setminus S$ denotes the graph obtained from $G$ by deleting all vertices in $S$ and their incident edges. $G[S]$ denotes the graph  induced by $S$ whose vertex set is $S$ and whose edge set consists of all edges of $G$ which have both ends in $S$. A set $M$ of disjoint edges of $G$ is called a {\sl matching} in $G$.  The {\sl matching number}, denoted by $\nu(G)$, is the maximum cardinality of a matching in $G$. We call a matching with $k$ edges a {\sl $k$-matching}, denoted by $M_k$.
For two vertex disjoint graphs $G$ and $H$, the {\sl union} of $G$ and $H$, denoted by $G\cup H$, is the graph with vertex set $V(G)\cup V(H)$ and edge set $E(G)\cup E(H)$. The {\sl join} of $G$ and $H$, denoted by $G\vee H$, is the graph obtained from $G\cup H$ by adding all edges between $V(G)$ and $V(H)$. $K_{r}$ denotes the complete graph of order $r$. Let $K_{r}(n_1,\cdots,n_r)$ be the complete $r$-partite graph with classes of sizes $n_1 ,\cdots , n_r$. If $\sum_{i=1}^{r}n_i$ $=n$ and $|n_i-n_j|\leq 1$ for any $1\leq i<j\leq r$, then $K_{r}(n_1,\cdots,n_r)$ is called an $r$-partite {\sl Tur\'{a}n graph}, denoted by $T_{r}(n)$.  Let $S_{k}$ be a star of $k$ edges. The {\sl center} of a star is the vertex of maximum degree in the star. A {\sl star forest} is the union of vertex disjoint stars. We use $S_{k_1,\cdots,k_q}$ to denote $\cup_{i=1}^{q}S_{k_i}$, where  $k_1\geq k_2\geq \cdots \geq k_q\geq 1$.

Given a family of graphs $\mathcal F$, a graph $G$ is called $\mathcal F$-{\sl free} if it does not contain any copy of $F$ in $\mathcal F$. If $\mathcal F=\{F\}$, we use $F$-free instead of $\mathcal F$-free. In 1941, Tur\'{a}n \cite{Turan1941} proved that if $G$ is $K_{r+1}$-free of order $n$, then $e(G)\leq e(T_{r}(n))$, with equality if and only if $G=T_{r}(n)$. The problem of determining the maximum number of edges among all $n$-vertex $\mathcal F$-free graphs is usually called the Tur\'{a}n type extremal problem. The {\sl Tur\'{a}n number}, denoted by $\mathrm{ex}(n,\mathcal F)$, is the maximum number of edges among all $n$-vertex $\mathcal F$-free graphs. An $n$-vertex $\mathcal F$-free graph with $\mathrm{ex}(n,\mathcal F)$ edges is called an {\sl extremal graph} for $\mathcal F$. Denote by $\mathrm{Ex}(n,\mathcal F)$ the set of all extremal graphs for $\mathcal F$.

Given a graph $G$ and an integer $p\geq 2$, the {\sl edge blow-up} of $G$, denoted by $G^{p+1}$, is obtained by replacing each edge of $G$ with a clique of order $p+1$, where the new vertices of the cliques are all distinct. Simonovits \cite{Simonovits1968} and Moon \cite{JWMoon1968} determined the extremal graphs for the edge blow-up of matchings. Chen, Gould, Pfender and Wei \cite{ChenGould2003} determined the Tur\'{a}n number for the edge blow-up of stars. There are other extremal results for the edge blow-up of paths, cycles and trees, see \cite{Glebov,Liu2013,WangHouLiuMa2021,NiKangShanZhu2020}. In 2022, Yuan \cite{Yuan2022} determined the range of the Tur\'{a}n numbers for edge blow-up of all bipartite graphs and the exact Tur\'{a}n numbers for edge blow-up of all non-bipartite graphs.

The {\sl spectral radius} of a graph $G$ is the spectral radius of its adjacent matrix $A(G)$, denoted by $\rho(G)$.
In 2010, Nikiforov \cite{Nikiforov2010extremalspectra} proposed the following spectral Tur\'{a}n type extremal problem: What is the maximum spectral radius among all $n$-vertex $\mathcal F$-free graphs? We denote by $\mathrm{Ex_{sp}}(n,F)$ the set of all graphs which attain the maximum spectral radius.
Nikiforov \cite{Nikiforov07} determined the spectral version of Tur\'{a}n's extremal result, and proved that if $G$ is a $K_{r+1}$-free graph on $n$ vertices, then $\rho (G)\le \rho (T_{r}(n))$, with equality if and only if $G=T_{r}(n)$.
Cioab\u{a}, Feng, Tait and Zhang \cite{CioabaFengTaitZhang} proved that the graph attaining the maximum spectral radius among all $n$-vertex $S_k^3$-free graphs is a member of $\mathrm{Ex}(n,S_k^3)$.
 In 2022, Zhai, Liu and Xue \cite{zhai2022} determined the unique spectral extremal graphs for $S_k^3$ when $n$ is sufficiently large. Let $H_{s,k}$ be the graph obtained from s triangles and $k$ odd cycles of lengths
at least $5$ by sharing a common vertex.
 Li and Peng \cite{Yongtao21} showed that $\mathrm{Ex_{sp}}(n, H_{s,k})\subseteq \mathrm{Ex}(n, H_{s,k})$ for $n$ large enough.
 Lin, Zhai and Zhao \cite{LZZ2022} proved that  $\mathrm{Ex_{sp}}(n,\Gamma_k)\subseteq \mathrm{Ex}(n,\Gamma_k)$ for sufficiently large $n$,
  where $\Gamma_k$  is the family of graphs consisting of $k$-edge-disjoint
triangles.
Desai, Kang, Li, Ni, Tait and Wang \cite{DKLNTW} generalized the result of \cite{CioabaFengTaitZhang} to $S_{k}^p$.
 Cioab\u{a}, Desai and Tait \cite{CDT21} raised the following conjecture: Let $F$ be any graph such that the graphs in $\mathrm{Ex}(n,F)$ are Tur\'{a}n graphs plus $O(1)$ edges. Then for sufficiently large $n$, a graph attaining the maximum spectral radius among all $F$-free graphs on $n$ vertices is a member of $\mathrm{Ex}(n,F)$.
 Wang, Kang and Xue \cite{Kang2023} proved the conjecture.
  Recently, Ni, Wang and Kang \cite{NiWangKang} proved that if $G$ has the maximum spectral radius among all $M^{p+1}_k$-free graphs on $n$ vertices, then $G\in \mathrm{Ex}(n,M^{p+1}_k)$.

 In this paper, we consider the spectral extremal problem for the edge blow-up of star forest $S^{p+1}_{k_1,\cdots,k_q}$.  The following result is obtained.
\begin{theorem}\label{main}
Let $p\geq 2$, $q\geq 1$, and $k_1\geq k_2\geq \cdots \geq k_q\geq 1$. For sufficiently large $n$, if $H$ has the maximum spectral radius among all $n$-vertex $S^{p+1}_{k_1,\cdots,k_q}$-free graphs, then  $$H\in \mathrm{Ex}(n,S^{p+1}_{k_1,\cdots,k_q}).$$
\end{theorem}

\section{Preliminaries}\label{sec2}

In this section, we give some notations and useful lemmas.
Let $G$ be a simple graph with matching number $\nu(G)$ and maximum degree $\Delta(G)$. For two given integers $\nu$ and $\Delta$, define $f(\nu, \Delta)
=\max\{e(G): \nu(G)\leq \nu, \Delta(G)\leq \Delta \}$. In 1976, Chv\'atal and Hanson \cite{Chvatal76} obtained the following result.

\begin{theorem}[\cite{Chvatal76}]\label{ffunction}
For integers $\nu \geq 1$ and $\Delta \geq 1$, we have
$$f(\nu, \Delta)= \Delta \nu +\left\lfloor\frac{\Delta}{2}\right\rfloor
 \left \lfloor \frac{\nu}{\lceil{\Delta}/{2}\rceil }\right \rfloor
 \leq \Delta \nu+\nu.$$
\end{theorem}
Denote by $\mathcal{E}_{\nu,\Delta}$  the set of the  graphs obtained the upper bound
in Theorem \ref{ffunction}.
Let $H(n, p, q)={K}_{q-1}\vee T_{p}(n-q+1)$ and $\mathcal{H}(n,p,q,k_q-1)$ be a set of  graphs obtained from $K_{q-1}\vee T_{p}(n-q+1)$ by putting a copy of  ${E}_{k_q-1, k_q-1}\in\mathcal{E}_{k_q-1, k_q-1}$ in one class of $T_{p}(n-q+1)$. Let $h(n,p,q)=e(H(n, p, q))$.
Using Theorem 2.3 in \cite{Yuan2022} on graph $G=S_{k_1,\cdots,k_q}$, we obtain that
\begin{equation*}
 \mathrm{ex}(n,S^{p+1}_{k_1,\cdots,k_q})\leq h(n,p,q)+f(k_q-1,k_q-1).
\end{equation*}
Combining with the fact that the graphs in $\mathcal{H}(n,p,q,k_q-1)$ are $S^{p+1}_{k_1,\cdots,k_q}$-free, we have the following result.

\begin{theorem}\label{extremal}
If $n$ is large enough, then we have
\begin{equation}\label{extremalnumber}
 \mathrm{ex}(n,S^{p+1}_{k_1,\cdots,k_q})= h(n,p,q)+f(k_q-1,k_q-1).
\end{equation}
Moreover, the graphs in  $\mathcal{H}(n,p,q,k_q-1)$ are  the only extremal graph for $S^{p+1}_{k_1,\cdots,k_q}$.
\end{theorem}


Chen, Gould, Pfender and Wei \cite{ChenGould2003} gave the following lemma.
\begin{lemma}[\cite{ChenGould2003}]\label{chen}
Suppose that $G$ is a $S_{k}^{r}$-free graph, and $V(G)=V_1\cup V_2\cup \cdots \cup V_{r-1}$ is a vertex partition of $G$. Let $E_{cr}(G)=\cup_{1\leq i<j\leq r-1}E(V_i,V_j)$ be the set of crossing edges of $G$.  For any $i\in [r-1]$ and $v\in V_i$, if
\begin{eqnarray*}
& &\sum_{j\neq i}\nu(G[V_j])\leq k-1\ \ \  \mbox{and} \ \ \ \Delta(G[V_i])\leq k-1,\\
& &d_{G[V_i]}(v)+\sum_{j\neq i}\nu(G[N(v)\cap V_j])\leq k-1,
\end{eqnarray*}
then $$\sum_{i=1}^{r-1}|E(G[V_i])|-\left(\sum_{1\leq i<j\leq r-1}|V_i||V_j|-|E_{cr}(G)|\right)\leq f(k-1,k-1).$$
\end{lemma}

The following lemma was given in \cite{CioabaFengTaitZhang}.
\begin{lemma}[\cite{CioabaFengTaitZhang}]
\label{intersect}
Let $V_1,\cdots,V_n$ be $n$ finite sets. Then
\[
|V_1 \cap \cdots \cap V_n| \geq \sum_{i=1}^{n}|V_i|-(n-1)|\bigcup_{i=1}^{n}V_i|
\]
\end{lemma}


By the Perron-Frobenius theorem, $\rho(G)$ is the largest eigenvalue of $A(G)$.
Let $\mathbf{x}=(x_1,\cdots,x_n)^{\mathrm{T}}$ be an eigenvector corresponding to $\rho(G)$. Then for any $ i\in [n]$ and $x_i\neq 0$,
\begin{equation}\label{eigenequation}
\rho(G)x_i=\sum_{ij\in E(G)}x_j.
\end{equation}
By the Rayleigh quotient, we have
\begin{equation}\label{Rayleigh}
\rho(G)=\max_{\mathbf{x}\in \mathbb{R}^{n}_{+}}\frac{\mathbf{x}^{\mathrm{T}}A(G)\mathbf{x}}{\mathbf{x}^{\mathrm{T}}\mathbf{x}}=\max_{\mathbf{x}\in \mathbb{R}^{n}_{+}}\frac{2\sum_{ij\in E(G)}x_ix_j}{\mathbf{x}^{\mathrm{T}}\mathbf{x}}.
\end{equation}

From the spectral version of the Stability Theorem given by Nikiforov \cite{Niki09JGT},  we get the following result.
\begin{lemma}[\cite{DKLNTW}] \label{stability}
Let $F$ be a graph with chromatic number $\chi (F)=r+1$. For every $\varepsilon >0$, there exist $\delta >0$ and $n_0$ such that if  $G$ is an $F$-free graph on $n\ge n_0$ vertices  with $\rho (G) \ge (1- \frac{1}{r} -\delta )n$, then $G$ can be obtained from $T_{r}(n)$ by adding and deleting at most $\varepsilon n^2$ edges.
\end{lemma}

\section{Proof of Theorem \ref{main}}
In this section, we suppose that $H$ has the maximum spectral radius among all $n$-vertex $S^{p+1}_{k_1,\cdots,k_q}$-free graphs. Then we will prove that $e(H)=\mathrm{ex}(n,S^{p+1}_{k_1,\cdots,k_q})$  for  sufficiently large $n$.
Clearly, $H$ is connected. Let $\rho(H)$ be the spectral radius of $H$. Then by the Perron-Frobenius theorem, there exist positive eigenvectors corresponding to $\rho(H)$.
Let $\mathbf{x}$ be a positive eigenvector corresponding to $\rho(H)$ with $\max\{x_i: i\in V(H)\}=1$. Without loss of generality, we assume that $x_z=1$.

\begin{lemma}\label{rho}
Let $H$ be an $n$-vertex $S^{p+1}_{k_1,\cdots,k_q}$-free graph with maximum spectral radius. Then
$$\rho(H)\geq \frac{p-1}{p}n+O(1).$$
\end{lemma}
\noindent{\bfseries Proof.}
Let $H'$ be an $n$-vertex $S^{p+1}_{k_1,\cdots,k_q}$-free graph with $e(H')=\mathrm{ex}(n,S^{p+1}_{k_1,\cdots,k_q})$.
By (\ref{extremalnumber}), we have
\begin{align*}
e(H')&\geq e(T_{p}(n-q+1))+(q-1)(n-q+1)+\binom{q-1}{2}\nonumber\\[2mm]
&\geq e(T_{p}(n))+\frac{q-1}{p}n-\frac{(q-1)(p+q-1)}{2p}-\frac{p}{8}.\label{exnumber}
\end{align*}
According to (\ref{Rayleigh}), for sufficiently $n$, we have
\begin{align*}
\rho(H)&\geq \rho(H')
\geq \frac{\mathbf{1}^{\mathrm{T}}A(H')\mathbf{1}}{\mathbf{1}^{\mathrm{T}}\mathbf{1}}
= \frac{2e(H')}{n}\\[2mm]
&\geq \frac{2}{n}\left(e(T_{p}(n))+\frac{q-1}{p}n-\frac{(q-1)(p+q-1)}{2p}-\frac{p}{8}\right)\\
&\geq\frac{p-1}{p}n+\frac{2(q-1)}{p}-\frac{1}{n}\left(\frac{(q-1)(p+q-1)}{p}+\frac{p}{2}\right)\\
&\geq\frac{p-1}{p}n+O(1).
\end{align*}
\qed

By Lemmas \ref{stability}, \ref{rho} and similar discussion as in
\cite{Kang2023}, we have the following lemmas.

\begin{lemma}\label{partition}
Let $H$ be an $n$-vertex $S^{p+1}_{k_1,\cdots,k_q}$-free graph with maximum spectral radius.
For every $\varepsilon >0$  
and sufficiently large $n$,
\[
e(H)\geq e(T_{p}(n))-\varepsilon n^2.
\]
Furthermore, there exists $\varepsilon_1 >0$ such that $H$ has a partition $V(H)=V_1 \cup \cdots \cup V_{p}$ with $\sum_{1\leq i<j\leq p}e(V_i,V_j)$ attaining the maximum,
\[
\sum_{i=1}^{p}e(V_i)\leq  \varepsilon n^2,
\]
and for any $i\in [p]$
$$\frac{n}{p}-\varepsilon_1 n< |V_i|< \frac{n}{p}+\varepsilon_1 n.$$
\end{lemma}


\begin{lemma}\label{W}
Suppose that $\varepsilon$ and $\theta$ are two sufficiently small constants with 
$\varepsilon\leq \theta^2$. Let
\[
W:=\cup_{i=1}^{p}\{v\in V_i: d_{V_i}(v)\geq 2 \theta n\}.
\]
Then $|W|\leq \theta n$.
\end{lemma}

\begin{lemma}\label{L}
Suppose that $\varepsilon_2$ is a sufficiently small constant with $\varepsilon < \varepsilon_2\ll \theta$. Let
\[
L:=\{v\in V(H): d(v)\leq (1-\frac{1}{p}-\varepsilon_2)n\}.
\]
Then $|L|\leq \varepsilon_3 n$, where $\varepsilon_3\ll \varepsilon_2$ is a sufficiently small constant satisfying $\varepsilon-\varepsilon_2 \varepsilon_3+\frac{p-1}{2p}\varepsilon_3^2<0$.
\end{lemma}
\noindent{\bfseries Proof.}
Suppose to the contrary that $|L|> \varepsilon_3 n$, then there exists $L'\subseteq L$ with $|L'|=\lfloor \varepsilon_3 n \rfloor$. Therefore,
\begin{eqnarray*}
e(H\setminus L')&\geq& e(H)-\sum_{v\in L'}d(v)\\[2mm]
&\geq & e(T_{p}(n))-\varepsilon n^2- \lfloor \varepsilon_3 n \rfloor (1-\frac{1}{p}-\varepsilon_2)n\\[2mm]
&> & \frac{p-1}{2p}(n-\lfloor \varepsilon_3 n \rfloor-q+1)^2+(q-1)(n-\lfloor \varepsilon_3 n \rfloor-q+1)+\binom{q-1}{2}+k_q^2\\[2mm]
&\geq  & e({K}_{q-1}\vee T_{p}(n-\lfloor \varepsilon_3 n \rfloor-q+1))+f(k_q-1,k_q-1)\\[2mm]
&= & \mathrm{ex}(n-\lfloor \varepsilon_3 n \rfloor,S^{p+1}_{k_1,\cdots,k_q}).
\end{eqnarray*}
Since $e(H\setminus L')>\mathrm{ex}(n-|L'|,S^{p+1}_{k_1,\cdots,k_q})$, $H\setminus L'$ contains a copy of $S^{p+1}_{k_1,\cdots,k_q}$. This contradicts the fact that $H$ is $S^{p+1}_{k_1,\cdots,k_q}$-free.
\qed

\begin{lemma}\label{Ttpp}
For an arbitrary subset $S\subseteq V_{i}\setminus (W\cup L)$, $i\in [p]$, there exists a subgraph $H'$ in $H$ such that $H'$ is a copy of
$T_{p}(p|S|)$ and $V(H')\cap V_{i}=S$.
\end{lemma}
\noindent{\bfseries Proof.}
For any $i\in [p]$ and any vertex $w\in V_i\setminus (W\cup L)$, we have $d(w)>(1-\frac{1}{p}-\varepsilon_2)n$ and $d_{V_i}(w)<2\theta n$. Then for any $j\in [p]$ and $j\neq i$, we have
\begin{eqnarray*}
 d_{V_j}(w)&\geq& d(w)-d_{V_i}(w)-(p-2)(\frac{n}{p}+\varepsilon_1 n)\\
&>& (1-\frac{1}{p}-\varepsilon_2)n- 2\theta n- (p-2)(\frac{n}{p}+\varepsilon_1 n)\\
&> &\frac{n}{p}-3\theta n-(p-2)\varepsilon_1 n.
\end{eqnarray*}
Take an arbitrary subset $S\subseteq V_{i}\setminus (W\cup L)$. Without loss of generality, let $i=1$ and $S=\{u_{1,1},\cdots, u_{1,t}\}$. We consider the common neighbors of $\{u_{1,1},\cdots, u_{1,t}\}$ in $V_2\setminus (W\cup L)$.
Combining with Lemma \ref{intersect}, we have
\begin{eqnarray*}
& &|(\cap_{j\in [t]} N_{V_{2}}(u_{1,j}))\setminus (W\cup L)|\\[2mm]
&\geq & \sum_{j=1}^{t}d_{V_{2}}(u_{1,j})-(t-1)|V_2|-|W|-|L|\\[2mm]
&> & t(\frac{n}{p}-3\theta n-(p-2)\varepsilon_1 n)-(t-1)(\frac{n}{p}+\varepsilon_1 n)-\theta n- \varepsilon_3 n\\[2mm]
&> & \frac{n}{p}-o(n)\\[2mm]
&> & t.
\end{eqnarray*}
So there exist vertices $u_{2,1},\cdots, u_{2,t}\in V_2\setminus (W\cup L)$ such that $\{u_{1,1},\cdots, u_{1,t}\} \cup \{u_{2,1},\cdots, u_{2,t}\}$ can induce a complete bipartite subgraph in $H$. For an integer $s$ with $2\leq s\leq p-1$, suppose that there are vertices $u_{s,1},\cdots,u_{s,t}\in V_{s}\setminus (W\cup L)$ such that $\{u_{1,1},\cdots, u_{1,t}\}~\cup$ $\{u_{2,1},\cdots,u_{2,t}\}~\cup$ $\cdots ~ \cup$ $\{u_{s,1},\cdots,u_{s,t}\}$  induce a complete $s$-partite subgraph in $H$. We next consider the common neighbors of the above $st$ vertices in $V_{s+1}\setminus (W\cup L)$.
By Lemma \ref{intersect}, we have
\begin{eqnarray*}
& &|(\cap_{i\in [s], j\in [t]} N_{V_{s+1}}(u_{i,j}))\setminus (W\cup L)|\\[2mm]
&\geq & \sum_{i=1}^{s}\sum_{j=1}^{t}d_{V_{s+1}}(u_{i,j})-(st-1)|V_{s+1}|-|W|-|L|\\[2mm]
&> & st(\frac{n}{p}-3\theta n-(p-2)\varepsilon_1 n)-(st-1)(\frac{n}{p}+\varepsilon_1 n)-\theta n- \varepsilon_3 n\\[2mm]
&> & \frac{n}{p}-o(n)\\[2mm]
&> & t.
\end{eqnarray*}
Then we can find vertices $u_{s+1,1},\cdots, u_{s+1,t}$ $\in V_{s+1}\setminus (W\cup L)$ such that $\{u_{1,1},\cdots, u_{1,t}\}~\cup$ $\{u_{2,1},\cdots,u_{2,t}\}~\cup$ $\cdots ~ \cup$ $\{u_{s+1,1},\cdots,u_{s+1,t}\}$  induce a complete $(s+1)$-partite subgraph in $H$. Therefore, for every $2\leq i\leq p$, there exist vertices $u_{p,1},\cdots,u_{p,t}$ $\in V_{i}\setminus (W\cup L)$ such that $\{u_{1,1},\cdots, u_{1,t}\}~\cup$ $\{u_{2,1},\cdots,u_{2,t}\}~\cup$ $\cdots ~\cup$ $\{u_{p,1},\cdots,u_{p,t}\}$  induce a complete $p$-partite subgraph $H'$ in $H$. The result follows.
\qed

\begin{lemma}\label{independent}
For each $i\in [p]$, $$\Delta(H[V_i\setminus (W\cup L)])< q(k_1 p+1),$$ $$\nu(H[V_i\setminus (W\cup L)])< \sum_{i\in [q]}k_i,$$ and there exists an independent set $I_i\subseteq V_i\setminus (W\cup L)$ such that $$|I_i|\geq |V_i\setminus (W\cup L)|-2(\sum_{i\in[q]}k_i-1).$$
\end{lemma}

\noindent{\bfseries Proof.}
We first prove that $\Delta(H[V_i\setminus (W\cup L)])<q(k_1 p+1)$ for any $i\in [p]$. Suppose to contrary that there exist $i_0\in [p]$ and  $u\in V_{i_0}\setminus (W\cup L)$ such that $d_{V_{i_0}\setminus (W\cup L)}(u)\geq q(k_1 p+1)$.
As $u\in V_{i_0}\setminus (W\cup L)$, there exists a vertex $w$ such  that $uw\notin E(H)$.
Let $H'$ be the graph with $V(H')=V(H)$ and $E(H')=E(H)\cup \{uw\}$. Then $H$ is a proper subgraph of $H'$. By the maximum of $\rho(H)$, $H'$ contains a copy of $S^{p+1}_{k_1,\cdots,k_q}$, say $H_1$. From the construction of $H'$, we see that $u\in V(H_1)$, and there is an $j_0\in [q]$ such that $H_1\setminus \{u\}$ contains  a copy of $\cup_{i\in [q]\setminus \{j_0\}}S^{p+1}_{k_i}$, say $H_2$. Obviously, $H_2\subseteq H$. Let $S=N_{V_{i_0}\setminus (W\cup L)}(u)$. 
Since $d_{V_{i_0}\setminus (W\cup L)}(u)\geq q(k_1 p+1)$, by Lemma \ref{Ttpp} there exists a copy of  $S^{p+1}_{k_{j_0}}$  with center $u$, denoted by $H_3$, such that $V(H_3)\cap V(H_2)=\emptyset$. Thus $H_2\cup H_3$ is a copy of $S^{p+1}_{k_1,\cdots,k_q}$ in $H$, which contradicts the fact that $H$ is $S^{p+1}_{k_1,\cdots,k_q}$-free.

Now we claim that for any $i\in [p]$, $\nu(H[V_i\setminus (W\cup L)])< \sum_{i\in [q]}k_i$. Otherwise, there exists $i_0\in [p]$ such that $H[V_{i_0}\setminus (W\cup L)]$ contains a matching $M$ with $\sum_{i\in [q]}k_i$ edges. By Lemma \ref{Ttpp}, there is a copy of $T_{p}(2|M|p)$ with $V(T_{p}(2|M|p))\cap V_{i_0}=V(M)$. Since $|V(M)|=2|M|=2\sum_{i\in [q]}k_i$, we can find a  copy of $S^{p+1}_{k_1,\cdots,k_q}$ in $H$,  a contradiction.

For every $i\in [p]$, let $M^i$ be a maximum matching of $H[V_i\setminus (W\cup L)]$, and $C^i$ be the set of vertices covered by $M^i$. Then $|C^i|\leq 2(\sum_{i\in [q]}k_i-1)$. By deleting all vertices of $C^i$, we can find  an independent set $I_i\subseteq V_i\setminus (W\cup L)$. Therefore, $|I_i|\geq |V_i\setminus (W\cup L)|-2(\sum_{i\in [q]}k_i-1)$.
\qed

\begin{lemma}\label{WminusL}
$|W\setminus L|\leq q-1$, and for any $u\in W\setminus L$, $d(u)=n-1$.
\end{lemma}
\noindent{\bfseries Proof.}
We first prove the following claim.

\noindent{\bfseries Claim A. }
For any $u\in W\setminus L$, $H$ contains a copy of  $S^{p+1}_{q(k_1 p+1)}$  with center $u$.

For any $u\in W\setminus L$, we have $d(u)>(1-\frac{1}{p}-\varepsilon_2)n$. Without loss of generality, we may assume that $u\in V_1$. Since $V(H)=V_1 \cup \cdots \cup V_{p}$ is the vertex partition that  maximizes the number of crossing edges of $H$,  $d_{V_1}(u)\leq \frac{1}{p}d(u)$. For any $2\leq j\leq p$,
\begin{eqnarray}
 d_{V_j}(u)&\geq& d(u)-d_{V_1}(u)-(p-2)(\frac{n}{p}+\varepsilon_1 n)\nonumber\\
&> &\frac{p-1}{p}(1-\frac{1}{p}-\varepsilon_2)n-(p-2)(\frac{n}{p}+\varepsilon_1 n)\nonumber\\
&>& \frac{n}{p^2}-\theta n-(p-2)\varepsilon_1 n.\label{eq5}
\end{eqnarray}
For any $i\in [p]$ and $v\in V_i \setminus (W\cup L)$, we have $d(v)>(1-\frac{1}{p}-\varepsilon_2)n$ and $d_{V_i}(v)<2\theta n$. Then for any $j\in [p]$ and $j\neq i$, we have
\begin{eqnarray}
d_{V_j}(v)&\geq& d(v)-d_{V_i}(v)-(p-2)(\frac{n}{p}+\varepsilon_1 n)\nonumber\\
&> &\frac{n}{p}-\varepsilon_2 n-2\theta n-(p-2)\varepsilon_1 n\nonumber\\
&> & \frac{n}{p}-3\theta n-(p-2)\varepsilon_1 n.\label{eq4}
\end{eqnarray}
By Lemmas \ref{W} and \ref{L}, we have
\begin{eqnarray*}
d_{V_1\setminus (W\cup L)}(u)&\geq& d_{V_1}(u)-|W\cup L|\\
&\geq & 2\theta n-\theta n-\varepsilon_2 n\\
&> & q(k_1 p+1).
\end{eqnarray*}
Let $u_{1,1},\cdots, u_{1,q(k_1 p+1)}$ be the neighbors of $u$ in $V_1\setminus (W\cup L)$. We consider the common neighbors of $u,u_{1,1},\cdots,u_{1,q(k_1 p+1)}$ in $V_2\setminus (W\cup L)$.
Combining with (\ref{eq5}), (\ref{eq4}) and  Lemma \ref{intersect},  we have
\begin{eqnarray*}
& &|N_{V_2}(u)\cap (\cap_{i\in[q(k_1 p+1)]} N_{V_2}(u_{1,i}))\setminus (W\cup L)|\\[2mm]
&\geq & d_{V_2}(u)+\sum_{i=1}^{q(k_1 p+1)}d_{V_2}(u_{1,i})-q(k_1 p+1)|V_2|-|W|-|L|\\[2mm]
&> &\frac{n}{p^2}-\theta n-(p-2)\varepsilon_1 n + q(k_1 p+1)(\frac{n}{p}-3\theta n-(p-2)\varepsilon_1 n)- q(k_1 p+1)(\frac{n}{p}+\varepsilon_1 n)-\theta n- \varepsilon_3 n\\[2mm]
&> & \frac{n}{p^2}-o(n)\\[2mm]
&> & q(k_1 p+1).
\end{eqnarray*}
Let $u_{2,1},\cdots, u_{2,q(k_1 p+1)}$ be the common neighbors of $\{u,u_{1,1},\cdots,u_{1,q(k_1 p+1)}\}$ in $V_2\setminus (W\cup L)$. For an integer $2\leq s\leq p-1$, suppose that $u_{s,1},\cdots,u_{s,q(k_1 p+1)}$ are the common neighbors of $\{u,u_{i,1},\cdots,u_{i,q(k_1 p+1)}: 1\leq i\leq s-1\}$ in $ V_{s}\setminus (W\cup L)$. We next consider the common neighbors of $\{u,u_{i,1},\cdots,u_{i,q(k_1 p+1)}: 1\leq i\leq s\}$ in $V_{s+1}\setminus (W\cup L)$.
By  (\ref{eq5}), (\ref{eq4}) and   Lemma \ref{intersect}, we have
\begin{eqnarray*}
& &|N_{V_{s+1}}(u)\cap (\cap_{i\in [s], j\in [q(k_1 p+1)]} N_{V_{s+1}}(u_{i,j}))\setminus (W\cup L)|\\[2mm]
&\geq & d_{V_{s+1}}(u)+\sum_{i=1}^{s}\sum_{j=1}^{q(k_1 p+1)}d_{V_{s+1}}(u_{i,j})-sq(k_1 p+1)|V_{s+1}|-|W|-|L|\\[2mm]
&> & \frac{n}{p^2}-\theta n-(p-2)\varepsilon_1 n+sq(k_1 p+1)(\frac{n}{p}-3\theta n-(p-2)\varepsilon_1 n)-sq(k_1 p+1)(\frac{n}{p}+\varepsilon_1 n)-\theta n- \varepsilon_3 n\\[2mm]
&> & \frac{n}{p^2}-o(n)\\[2mm]
&> &q(k_1 p+1).
\end{eqnarray*}
Let $u_{s+1,1},\cdots, u_{s+1,q(k_1 p+1)}$ be the common neighbors of $\{u,u_{i,1},\cdots,u_{i,q(k_1 p+1)}: 1\leq i\leq s\}$ in $V_{s+1}\setminus (W\cup L)$. Therefore, for every $i\in [p]$, there exist $u_{i,1},\cdots,u_{i,q(k_1 p+1)}\in V_{i}\setminus (W\cup L)$ such that $\{u_{1,1},\cdots,u_{1,q(k_1 p+1)}\}$ $\cup ~\{u_{2,1},\cdots,u_{2,q(k_1 p+1)}\}$ $\cup ~\cdots ~\cup$ $\{u_{p,1},\cdots,u_{p,q(k_1 p+1)}\}$ form a complete $p$-partite subgraph of $H$, and $u$ is adjacent to all the above $pq(k_1 p+1)$ vertices. Hence we can find a copy of  $S^{p+1}_{q(k_1 p+1)}$  with center $u$.

Now we  prove that $|W\setminus L|\leq q-1$. Suppose to the contrary that $|W\setminus L|\geq q$. By Claim  A, for any $u\in W\setminus L$, we can find a copy of $S^{p+1}_{q(k_1 p+1)}$  with center $u$ in $H$.  Therefore, we can find $q$ disjoint copies of $S^{p+1}_{k_1}$  in $H$. This contradicts the fact that $H$ is $S^{p+1}_{k_1,\cdots,k_q}$-free.

Next we prove that for any $u\in W\setminus L$, $d(u)=n-1$. Suppose to the contrary that there exists $u\in W\setminus L$ such that $d(u)<n-1$. Assume that $v$ is a vertex of $H$ such that $uv\notin E(H)$. Let $H'$ be the graph with $V(H')=V(H)$ and $E(H')=E(H)\cup \{uv\}$. We claim that $H'$ is $S^{p+1}_{k_1,\cdots,k_q}$-free. Otherwise, $H'$ contains a copy of   $S^{p+1}_{k_1,\cdots,k_q}$, say $H_1$, then $uv\in E(H_1)$. Therefore, $H_1$ contains a copy of $S^{p+1}_{k_{i_0}}$, say $H_2$, such that $uv\in E(H_2)$, where $i_0\in [q]$.
Let $H_3=H_1\setminus H_2$. Then $H_3$ is a copy of  $\cup_{i\in[q]\setminus \{i_0\}}S^{p+1}_{k_i}$ in $H$. Since $u\in  W\setminus L$, we can find a copy of $S^{p+1}_{q(k_1 p+1)}$  with center $u$ in $H$ by Claim A. Thus there is a copy of  $S^{p+1}_{k_{i_0}}$  with center $u$ in $H$, denoted by $H_4$, such that $V(H_4)\cap V(H_3)=\emptyset$. Then $H_3\cup H_4$ is a  copy of $S^{p+1}_{k_1,\cdots,k_q}$ in $H$, and this is a contradiction. Therefore, $H'$ is $S^{p+1}_{k_1,\cdots,k_q}$-free. By the construction of $H'$, we have $\rho(H')>\rho(H)$, which contradicts the assumption that $H$ has the maximum spectral radius among all $n$-vertex $S^{p+1}_{k_1,\cdots,k_q}$-free graphs.
\qed

\begin{lemma}\label{xv0}
Let $x_{v_0}=\max\{x_v : v\in V(H)\setminus W\}$. Then $x_{v_0}>1-\frac{2}{p}$ and  $v_0\notin L$.

\end{lemma}
\noindent{\bfseries Proof.}
Recall that $x_z=\max\{x_v : v\in V(H)\}=1$, then
\[
\rho(H)=\rho(H)x_z\leq |W|+(n-|W|)x_{v_0}.
\]
By Lemmas \ref{L} and \ref{WminusL}, we have
\begin{eqnarray}\label{Wnumber}
|W|=|W\cap L|+|W\setminus L|\leq |L|+q-1\leq \varepsilon_3 n+q-1. \label{5.1}
\end{eqnarray}
Combining with Lemma \ref{rho}, we have
\begin{eqnarray}
x_{v_0}\geq \frac{\rho(H)-|W|}{n-|W|}
\geq \frac{\rho(H)-|W|}{n}
\geq 1-\frac{1}{p}- \varepsilon_3-\frac{O(1)}{n}
>1-\frac{2}{p}.\label{6.1}
\end{eqnarray}
Therefore, we have
\begin{eqnarray}
\rho (H) x_{v_0}&=&\sum_{v\sim v_0} x_v
= \sum_{\substack{v\in W, \\ v\sim v_0}} x_v
+ \sum_{\substack{v\notin W,\\ v\sim v_0} } x_v \nonumber\\[2mm]
&\leq & |W|+(d(v_0)-|W|)x_{v_0}. \label{7.1}
\end{eqnarray}

By (\ref{5.1}), (\ref{6.1}), (\ref{7.1}) and Lemma \ref{rho}, we have
\begin{eqnarray*}
d(v_0)&\geq& \rho(H)+|W|-\frac{|W|}{x_{v_0}}\\
&\geq& \rho(H)-\frac{2|W|}{p-2}\\
&\geq& \frac{p-1}{p}n+O(1)-\frac{2\varepsilon_3 n}{p-2}-\frac{2(q-1)}{p-2}\\
&>& (1-\frac{1}{p}-\varepsilon_2)n,
\end{eqnarray*}
where the last inequality  holds as $\varepsilon_3\ll \varepsilon_2$. Thus $v_0\notin L$.
\qed

\begin{lemma}\label{Lemptyset}
$L=\emptyset.$
\end{lemma}
\noindent{\bfseries Proof.}
Suppose to the contrary that there is a vertex $u_0\in L$, then $d(u_0)\leq (1-\frac{1}{p}-\varepsilon_2)n$. 
Let $H'$ be the graph with $V(H')=V(H)$ and $E(H')=E(H\setminus \{u_0\}) \cup \{wu_0: w\in \cup_{i=2}^{p}I_i\}$. It is obvious that $H'$ is $S^{p+1}_{k_1,\cdots,k_q}$-free. By Lemma \ref{xv0}, $v_0\in V(H)\setminus (W\cup L)$. Without loss of generality, we assume that $v_0\in V_1\setminus (W\cup L)$. Combining with Lemma \ref{independent}, we have
\begin{eqnarray*}
 \rho (H) x_{v_0}
 &=& \sum_{\substack{v\in W\cup L,\\ v\sim v_0}} x_v  + \sum_{\substack{v\in V_1\setminus (W\cup L),\\ v\sim v_0}} x_v
 + \sum_{\substack{v\in (\cup_{i=2}^{p}V_i)\setminus (W\cup L), \\ v\sim v_0} } x_v\\[2mm]
& < & |W|+|L|x_{v_0}+q(k_1 p+1)x_{v_0}+ \sum_{\substack{v\in (\cup_{i=2}^{p}V_i\setminus I_i)\setminus (W\cup L), \\ v\sim v_0} } x_v + \sum_{\substack{v\in \cup_{i=2}^{p}I_i, \\ v\sim v_0} } x_v\\[2mm]
& \leq & |W|+|L|x_{v_0}+q(k_1 p+1)x_{v_0}+ 2(\sum_{i\in [q]}k_i-1)(p-1)x_{v_0}+ \sum_{\substack{v\in \cup_{i=2}^{p}I_i } } x_v,
\end{eqnarray*}
which implies that
\begin{equation*}
\sum_{\substack{v\in \cup_{i=2}^{p}I_i } } x_v\geq (\rho (H)-|L|-q(k_1 p+1)-2(\sum_{i\in [q]}k_i-1)(p-1)) x_{v_0}-|W|.
\end{equation*}
Combining with Lemmas \ref{rho}, \ref{L} and \ref{xv0}, 
we have
\begin{eqnarray*}
&&\rho(H') - \rho(H)\\[2mm]
&\geq& \frac{\mathbf{x}^T\left(A(H')-A(H)\right)\mathbf{x}}{\mathbf{x}^T\mathbf{x}}\\[2mm]
& =& \frac{2x_{u_0}}{\mathbf{x}^T\mathbf{x}}\left(\sum_{\substack{w\in \cup_{i=2}^{p}I_i}} x_w - \sum_{u\sim u_0} x_u\right) \\[2mm]
& \geq& \frac{2x_{u_0}}{\mathbf{x}^T\mathbf{x}}\Bigl((\rho (H)-|L|-q(k_1 p+1)-2(\sum_{i\in [q]}k_i-1)(p-1)) x_{v_0}-2|W|-(d(u_0)-|W|)x_{v_0} \Bigr)\\[2mm]
& = &\frac{2x_{u_0}}{\mathbf{x}^T\mathbf{x}}\Bigl((\rho (H)-|L|-q(k_1 p+1)-2(\sum_{i\in [q]}k_i-1)(p-1)-d(u_0)+|W|) x_{v_0}-2|W| \Bigr)\\[2mm]
&\geq& \frac{2x_{u_0}}{\mathbf{x}^T\mathbf{x}}\Bigl( \frac{p-2}{p}(\varepsilon_2 n- \varepsilon_3 n- O(1))-2(\varepsilon_3 n+q-1) \Bigr)>0,
\end{eqnarray*}
where the last inequality holds since $\varepsilon_3 \ll \varepsilon_2$. This contradicts the fact that $H$ has the largest spectral radius over all $n$-vertex $S^{p+1}_{k_1,\cdots,k_q}$-free graphs. So $L$ must be empty.
\qed

\begin{lemma}\label{eigenvector}
For any $v\in V(H)$, $x_v\geq x_{v_0}-\frac{100k_1pq}{n}$.
\end{lemma}
\noindent{\bfseries Proof.}
Suppose to the contrary that there exists $u\in V(H)$ such that $x_u< x_{v_0}-\frac{100k_1pq}{n}$. Let $H'$ be the graph with $V(H')=V(H)$ and $E(H')=E(H\setminus \{u\})\cup \{uw: w\in \cup_{i=2}^{p}I_i\}$. It is obvious that $H'$ is $S^{p+1}_{k_1,\cdots,k_q}$-free.
Since $L=\emptyset$, then $|W|=|W\setminus L|\leq q-1$ by Lemma \ref{WminusL}. Recall that $x_{v_0}=\max\{x_v : v\in V(H)\setminus W\}$, then
 \begin{eqnarray*}
 \rho (H) x_{v_0}
 &=& \sum_{\substack{v\in W,\\ v\sim v_0}} x_v  + \sum_{\substack{v\in V_1\setminus W, \\ v\sim v_0}} x_v
 + \sum_{\substack{v\in (\cup_{i=2}^{p}V_i)\setminus W, \\ v\sim v_0} } x_v\\[2mm]
& < & |W|+q(k_1 p+1)x_{v_0}+  \sum_{\substack{v\in (\cup_{i=2}^{p}V_i\setminus I_i)\setminus W,\\ v\sim v_0} } x_v + \sum_{\substack{v\in \cup_{i=2}^{p}I_i, \\ v\sim v_0} } x_v\\[2mm]
& \leq & q-1+q(k_1 p+1)x_{v_0}+ 2(\sum_{i\in [q]}k_i-1)(p-1)x_{v_0}+ \sum_{\substack{v\in \cup_{i=2}^{p}I_i } } x_v,
\end{eqnarray*}
which implies that
\begin{equation*}
\sum_{\substack{v\in \cup_{i=2}^{p}I_i } } x_v\geq (\rho (H)-q(k_1 p+1)-2(\sum_{i\in [q]}k_i-1)(p-1)) x_{v_0}-(q-1).
\end{equation*}
 Therefore, we have
\begin{eqnarray*}
&&\rho(H') - \rho(H)\\
& \geq& \frac{2x_{u}}{\mathbf{x}^T\mathbf{x}}\Bigl((\rho (H)-q(k_1 p+1)-2(\sum_{i\in [q]}k_i-1)(p-1)) x_{v_0}-(q-1)- \rho(H)x_u \Bigr)\\[2mm]
&> &\frac{2x_{u}}{\mathbf{x}^T\mathbf{x}}\Bigl((\rho (H)-q(k_1 p+1)-2(\sum_{i\in [q]}k_i-1)(p-1)) x_{v_0} -(q-1)- \rho(H)(x_{v_0}-\frac{100k_1pq}{n}) \Bigr)\\[2mm]
&> &\frac{2x_{u}}{\mathbf{x}^T\mathbf{x}}\Bigl( \frac{(p-1)n}{p}\frac{100k_1pq}{n}- q(k_1 p+1)-2(\sum_{i\in [q]}k_i-1)(p-1)-(q-1)\Bigr)>0.
\end{eqnarray*}
This contradicts the fact that $H$ has the largest spectral radius over all $n$-vertex $S^{p+1}_{k_1,\cdots,k_q}$-free graphs.
\qed

\begin{lemma}\label{Wk-1}
$|W|=q-1$.
\end{lemma}
\noindent{\bfseries Proof.}
Let $|W|=s$. Then $s\leq q-1$ by Lemmas \ref{WminusL} and \ref{Lemptyset}. Suppose to the contrary that $s<q-1$.
By Lemmas \ref{independent} and \ref{Lemptyset}, for any $i\in [p]$, $\nu(H[V_i\setminus W])<\sum_{i\in [q]}k_i$ and $\Delta(H[V_i\setminus W])<q(k_1 p+1)$. Thus $\nu(\cup_{i=1}^{p}H[V_i\setminus W])<p\sum_{i\in [q]}k_i$ and $\Delta(\cup_{i=1}^{p}H[V_i\setminus W])<q(k_1 p+1)$. Combining with Lemma \ref{ffunction}, we have
\begin{align*}
e(\cup_{i=1}^{p}H[V_i\setminus W])&\leq f(\nu(\cup_{i=1}^{p}H[V_i\setminus W]),\Delta(\cup_{i=1}^{p}H[V_i\setminus W]))\\
&< f(p\sum_{i\in [q]}k_i,q(k_1 p+1))\\
&< q(k_1 p+1)(k_1pq+1).
\end{align*}
Take $S\subseteq V_1\setminus W$ with $|S|=q-s-1$. Let $H'$ be the graph with $V(H')=V(H)$ and $E(H')=E(H)\setminus \{uv: uv\in \cup_{i=1}^{p}E(H[V_i\setminus W])\}\cup \{uv: u\in S,v\in (V_1\setminus W)\setminus S\}$. It is obvious that $H'$ is $S^{p+1}_{k_1,\cdots,k_q}$-free. Therefore,
\begin{align*}
\rho(H') - \rho(H)
&\geq \frac{\mathbf{x}^T\left(A(H')-A(H)\right)\mathbf{x}}{\mathbf{x}^T\mathbf{x}}\\[2mm]
&=\frac{2}{\mathbf{x}^T\mathbf{x}}\left(\sum_{ij\in E(H')}x_ix_j-\sum_{ij\in E(H)}x_ix_j\right) \\[2mm]
& \geq \frac{2}{\mathbf{x}^T\mathbf{x}}\left(|S|(|V_1|-|W|-|S|)(x_{v_0}-\frac{100k_1pq}{n})^2- q(k_1 p+1)(k_1pq+1)\right) \\[2mm]
& \geq \frac{2}{\mathbf{x}^T\mathbf{x}}\Bigl( (q-s-1)(\frac{n}{p}-o(n))(\frac{p-2}{p}-o(1))^2-q(k_1 p+1)(k_1pq+1)\Bigr)\\[2mm]
&>0.
\end{align*}
This contradicts the fact that $H$ has the largest spectral radius over all $n$-vertex $S^{p+1}_{k_1,\cdots,k_q}$-free graphs.
Therefore, $|W|=s=q-1$.
\qed




\begin{lemma}\label{Bi}
For any $i\in [p]$, let $B_i=\{u\in V_i\setminus W: d_{V_i\setminus W}(u)\geq 1\}$ and $C_i=(V_i\setminus W)\setminus B_i$. Then $|B_i|< 2k_q^2$, and
 $u$ is adjacent to all vertices of $V(H)\setminus V_i$ for every vertex $u\in C_i$.
\end{lemma}
\noindent{\bfseries Proof.}
We first prove the following claim.

\noindent{\bfseries Claim B.}
For any $i\in [p]$, $H[V_i\setminus W]$ is both $S_{k_q}$-free and $M_{k_q}$-free.

Suppose to the contrary that there is an $i_0\in [p]$ such that $H[V_{i_0}\setminus W]$ contains a $S_{k_q}$ or an $M_{k_q}$ as a subgraph. Then $H$ contains  a copy of $S^{p+1}_{k_q}$, say $H_0$, by Lemmas \ref{Ttpp} and \ref{Lemptyset}. Since $|W|=q-1$ and each vertex of $W$ has degree $n-1$, $H$ contains a copy of $S^{p+1}_{k_1,\cdots,k_{q-1}}$, say $H_0'$, such that $V(H_0)\cap V(H_0')=\emptyset$. Thus $H_0\cup H_0'$ is a copy of  $S^{p+1}_{k_1,\cdots,k_{q}}$ in  $H$,  a contradiction.

By Claim B, we have
 $e(H[V_i\setminus W])\leq f(k_q-1,k_q-1)<k_q^2$. Therefore,
\[
|B_i|\leq  \sum_{u\in B_i}d_{V_i}(u)=\sum_{u\in V_i\setminus W}d_{V_i}(u)=2e(H[V_i\setminus W])<2k_q^2.
\]

Now we prove that each vertex of $C_i$ is adjacent to all vertices of $V\setminus V_i$ for any $i\in [p]$. Since each vertex of $W$ has degree $n-1$, we only need to prove that each vertex of $C_i$ is adjacent to all vertices of $(V(H)\setminus W)\setminus V_i$ for any $i\in [p]$. Suppose to the contrary that there exist $i_0\in [p]$ and $v\in C_{i_0}$ such that there is a vertex $w_{1,1}\notin V_{i_0}\cup W$ and $vw_{1,1}\notin E(H)$. Without loss of generality, we may assume that $i_0=1$. Let $H'$ be the graph with $V(H')=V(H)$ and $E(H')=E(H)\cup \{vw_{1,1}\}$. We claim that $H'$ is $S^{p+1}_{k_1,\cdots,k_{q}}$-free. Otherwise, $H'$ contains a  copy of $S^{p+1}_{k_1,\cdots,k_{q}}$, say $H_1$. Then $H_1$ has a copy of  $S^{p+1}_{k_{j_0}}$, say $H_2$, such that $vw_{1,1}\in E(H_2)$, where $j_0\in [q]$. Without loss of generality, let $j_0=q$. Let $H_3=H_1\setminus H_2$.
We first claim that $V(H_2)\cap W=\emptyset$. Otherwise, by Claim A in Lemma \ref{WminusL}, we can find a copy of $S^{p+1}_{k_{q}}$ in $H$, say $H_4$, such that $V(H_4)\cap V(H_3)=\emptyset$ and $H_3\cup H_4$ is a copy of $S^{p+1}_{k_1,\cdots,k_{q}}$ in $H$, a contradiction.

We may assume that $v$ is the center of $H_2$ (the case that $v$ is not the center of $H_2$ can be proved similarly). As $v$ is the center of $H_2$ and $V(H_2)\cap W=\emptyset$, there exist vertices $w_{1,1}, w_{1,2},\cdots, w_{1,p}, w_{2,1},\cdots, w_{2,p},\cdots,$ $w_{k_q,1},\cdots,w_{k_q,p}\notin V_1\setminus W$ such that for any $i\in [k_q]$, the vertex set $\{w_{i,1},w_{i,2},\cdots, w_{i,p}\}$ induces a $K_{p}$  in $H$.  For any $i\in [k_q]$ and $j\in [p]$, there is $2\leq s\leq p$ such that $w_{i,j}\in V_s$, then
\begin{align*}
d_{V_1}(w_{i,j})&=d(w_{i,j})-d_{V_s}(w_{i,j})-(p-2)(\frac{n}{p}+ \epsilon_1 n)\\
&\geq (\frac{p-1}{p}-\epsilon_2)n-(q-1+k_q-1)-(p-2)(\frac{n}{p}+ \epsilon_1 n)\\
&> \frac{n}{p}- \epsilon_2n -p\epsilon_1 n,
\end{align*}
where the last second inequality holds as $L=\emptyset$, $|W|=q-1$ and $H[V_s\setminus W]$ is $S_{k_q}$-free.
Using Lemma \ref{intersect}, we have
\begin{eqnarray*}
& & |\cap_{i\in [k_q],j\in [p]}N_{V_1}(w_{i,j})\setminus (W\cup B_1)|\\[2mm]
&\geq &\sum_{i\in [k_q],j\in [p]}d_{V_1}(w_{i,j})-(k_qp-1)|V_1|- |W|- |B_1|\\[2mm]
&> & k_qp(\frac{n}{p}- \epsilon_2n -p\epsilon_1 n)-(k_qp-1)(\frac{n}{p}+ \epsilon_1 n)-(q-1)-2k_q^2\\[2mm]
&\geq & \frac{n}{p}-o(n)\\[2mm]
& >& q(k_1p+1).
\end{eqnarray*}
Then there exists $v'\in C_1$ and $v'\notin V(H_1)$ such that $(H_2\setminus \{v\})\cup \{v'\}$ is a copy of   $S^{p+1}_{k_q}$, say $H_5$, such that $V(H_3)\cap V(H_5)=\emptyset$. Thus $H_3\cup H_5$ is a copy of $S^{p+1}_{k_1,\cdots,k_{q}}$ in $H$,  a contradiction. Therefore, $H'$ is $S^{p+1}_{k_1,\cdots,k_q}$-free. From the construction of $H'$, we see that $\rho(H')>\rho(H)$, which contradicts the assumption that $H$ has the maximum spectral radius among all $n$-vertex $S^{p+1}_{k_1,\cdots,k_q}$-free graphs.
\qed


\begin{lemma}\label{balance}
For any $1\leq i<j\leq p$,  $\left||V_i\setminus W|-|V_j\setminus W|\right|\leq 1$.
\end{lemma}

\noindent{\bfseries Proof.}
For any $i\in [p]$, let $|V_i\setminus W|=n_i$. Suppose that $n_1\geq n_2\geq \cdots \geq n_{p}$.  We prove the assertion by contradiction.  Assume that there exist $i_0, j_0$ with $1\leq i_0 < j_0\leq p$ such that $n_{i_0}-n_{j_0}\geq 2$. Let $F=K_p(n_1,n_2,\cdots,n_{p})$ 
 and $F'=K_{p}(n_1,\cdots, n_{i_0}-1,\cdots,n_{j_0}+1,\cdots,n_{p}).$ Assume that $F'\cong K_{p}(n'_1,n'_2,\cdots,n'_{p})$, where $n'_1\geq n'_2\geq \cdots \geq n'_{p}$.

\noindent{\bfseries Claim 1.} There exists a constant $c_1>0$ such that $$\rho(K_{q-1}\vee T_{p}(n-q+1))-\rho(K_{q-1}\vee F)\geq \frac{c_1}{n}.$$

Let $\mathbf{y}$ be a positive  eigenvector of $K_{q-1}\vee F$ corresponding to $\rho(K_{q-1}\vee F)$. By the symmetry we may assume $\mathbf{y}=(\underbrace{y_1,\cdots,y_1}_{n_1},\underbrace{y_2,\cdots,y_2}_{n_2},\cdots, \underbrace{y_p,\cdots,y_p}_{n_p},\underbrace{y_{p+1},\cdots,y_{p+1}}_{q-1})^{\mathrm{T}}$.
Thus by  (\ref{eigenequation}), we have
\begin{eqnarray}
\rho(K_{q-1}\vee F) y_i=\sum_{j=1}^{p}n_jy_j-n_iy_i+(q-1)y_{p+1}, \text{ for any } i\in [p], \label{114}
\end{eqnarray}
and \begin{eqnarray}
\rho(K_{q-1}\vee F) y_{p+1}=\sum_{j=1}^{p}n_jy_j+(q-2)y_{p+1}.\label{115}
\end{eqnarray}
Combining  (\ref{114}) and (\ref{115}), we have $y_i=\frac{\rho(K_{q-1}\vee F)+1}{\rho(K_{q-1}\vee F)+n_i}y_{p+1}$ for any $i\in [p]$, which implies that $y_{p+1}=\max\{y_v : v\in V(K_{q-1}\vee F)\}$. Without loss of generality, we assume that $y_{p+1}=1$, then
\begin{eqnarray}
y_i=\frac{\rho(K_{q-1}\vee F)+1}{\rho(K_{q-1}\vee F)+n_i}, \text{ for any } i\in [p]. \label{116}
\end{eqnarray}
Let $u_{i_0}\in V_{i_0}\setminus W$ be a fixed vertex. Then $K_{q-1}\vee F'$ can be obtained from $K_{q-1}\vee F$ by deleting all edges between $u_{i_0}$ and $V_{j_0}\setminus W$, and adding all edges between $u_{i_0}$ and $V_{i_0}\setminus (W\cup\{u_{i_0}\})$.
According to  (\ref{Rayleigh}), we deduce that
\begin{eqnarray*}
&&\rho(K_{q-1}\vee F')-\rho(K_{q-1}\vee F)\\[2mm]
&\geq& \frac{\mathbf{y}^{\mathrm{T}}(A(K_{q-1}\vee F')-A(K_{q-1}\vee F))\mathbf{y}}{\mathbf{y}^{\mathrm{T}}\mathbf{y}}\\[2mm]
&=& \frac{2}{\mathbf{y}^{\mathrm{T}}\mathbf{y}}\left((n_{i_0}-1)y_{i_0}^2-n_{j_0}y_{i_0}y_{j_0}\right)\\[2mm]
&=& \frac{2y_{i_0}}{\mathbf{y}^{\mathrm{T}}\mathbf{y}}\left((n_{i_0}-1)\frac{\rho(K_{q-1}\vee F)+1}{\rho(K_{q-1}\vee F)+n_{i_0}}-n_{j_0}\frac{\rho(K_{q-1}\vee F)+1}{\rho(K_{q-1}\vee F)+n_{j_0}}\right)\\[2mm]
&=& \frac{2y_{i_0}}{\mathbf{y}^{\mathrm{T}}\mathbf{y}}\frac{(\rho(K_{q-1}\vee F)+1)(n_{i_0}\rho(K_{q-1}\vee F)-n_{j_0}\rho(K_{q-1}\vee F)-\rho(K_{q-1}\vee F)-n_{j_0})}{(\rho(K_{q-1}\vee F)+n_{i_0})(\rho(K_{q-1}\vee F)+n_{j_0})}\\[2mm]
&\geq& \frac{2y_{i_0}}{\mathbf{y}^{\mathrm{T}}\mathbf{y}}\frac{(\rho(K_{q-1}\vee F)+1)(\rho(K_{q-1}\vee F)-n_{j_0})}{(\rho(K_{q-1}\vee F)+n_{i_0})(\rho(K_{q-1}\vee F)+n_{j_0})},
\end{eqnarray*}
where the last inequality holds as $n_{i_0}-n_{j_0}\geq 2$. Since $\delta(K_{q-1}\vee F)\leq \rho(K_{q-1}\vee F)\leq \Delta(K_{q-1}\vee F)$, in view of  the construction of $K_{q-1}\vee F$, we see that $\rho(K_{q-1}\vee F)=\Theta(n)$.
Combining  (\ref{116}) and the fact $\mathbf{y}^{\mathrm{T}}\mathbf{y}\leq n$, there exists a constant $c_1>0$ such that $\rho(K_{q-1}\vee F')-\rho(K_{q-1}\vee F)\geq \frac{c_1}{n}$, and this implies that $\rho(K_{q-1}\vee T_{p}(n-q+1)\geq \rho(K_{q-1}\vee F')$. Therefore, $\rho(K_{q-1}\vee T_{p}(n-q+1))-\rho(K_{q-1}\vee F)\geq\rho(K_{q-1}\vee F')-\rho(K_{q-1}\vee F)\geq \frac{c_1}{n}$.

\noindent{\bfseries Claim 2.} There exists a constant $c_2>0$ such that $$\rho(H)-\rho(K_{q-1}\vee T_{p}(n-q+1))\geq \frac{2f(k_q-1,k_q-1)}{n}\left(1-\frac{c_2}{n}\right).$$

Suppose that $\mathbf{z}$ is an eigenvector corresponding to $\rho(K_{q-1}\vee T_{p}(n-q+1))$ with $\max\{z_v : v\in V(K_{q-1}\vee T_{p}(n-q+1))\}=1$. Let  $n_1=\lfloor\frac{n-q+1}{p}\rfloor$, $n_2=\lceil\frac{n-q+1}{p}\rceil$, and $a=n-q+1-pn_1$. By the symmetry we may assume $\mathbf{z}=(\underbrace{z_1,\cdots,z_1}_{(p-a)n_1},\underbrace{z_2,\cdots,z_2}_{an_2}, \underbrace{z_{3},\cdots,z_{3}}_{q-1})^{\mathrm{T}}$. By the similar
discussion as in Claim 1,
 we have $z_3=1$,
\[
z_1= \frac{\rho(K_{q-1}\vee T_{p}(n-q+1))+1}{\rho(K_{q-1}\vee T_{p}(n-q+1))+n_1}>
z_2= \frac{\rho(K_{q-1}\vee T_{p}(n-q+1))+1}{\rho(K_{q-1}\vee T_{p}(n-q+1))+n_2}> 1-\frac{1}{p}.
\]
Let $H'$ be a graph obtained from $K_{q-1}\vee T_{p}(n-q+1)$ by embedding a copy of $E_{k_q-1,k_q-1}\in \mathcal{E}_{k_q-1,k_q-1}$ in one class of $T_{p}(n-q+1)$. It is obvious that $H'$ is $S^{p+1}_{k_1,\cdots,k_{q}}$-free. Therefore
\begin{align*}
\rho(H)-\rho(K_{q-1}\vee T_{p}(n-q+1))&\geq \rho(H')-\rho(K_{q-1}\vee T_{p}(n-q+1))\nonumber\\[2mm]
&\geq \frac{\mathbf{z}^\mathrm{T}(A(H')-A(K_{q-1}\vee T_{p}(n-q+1)))\mathbf{z}}{\mathbf{z}^{\mathrm{T}}\mathbf{z}}\nonumber\\[2mm]
&\geq \frac{2\sum_{ij\in E(E_{k_q-1,k_q-1})}z_{i}z_{j}}{\mathbf{z}^{\mathrm{T}}\mathbf{z}}\nonumber\\[2mm]
&\geq  \frac{2f(k_q-1,k_q-1)z_2^2}{\mathbf{z}^{\mathrm{T}}\mathbf{z}}.
\end{align*}
Since $n_2-n_1\leq 1$, $an_2+(p-a)n_1=n-q+1$ and $\rho(K_{q-1}\vee T_{p}(n-q+1))\geq \delta(K_{q-1}\vee T_{p}(n-q+1))=n-n_2$, we have
\begin{align*}
\frac{\mathbf{z}^{\mathrm{T}}\mathbf{z}}{z_2^2}&\leq \frac{q-1}{z_2^2}+ \frac{an_2z_1^2}{z_2^2}+ (p-a)n_1\\[2mm]
&= \frac{q-1}{z_2^2}+an_2\left(\frac{\rho(K_{q-1}\vee T_{p}(n-q+1))+n_2}{\rho(K_{q-1}\vee T_{p}(n-q+1))+n_1}\right)^2+(p-a)n_1\\[2mm]
&= \frac{q-1}{z_2^2}+an_2+(p-a)n_1+an_2\left(\left(1+\frac{1}{\rho(K_{q-1}\vee T_{p}(n-q+1))+n_1}\right)^2-1\right)\\[2mm]
&\leq  n+ (q-1)\left(\frac{1}{z_2^2}-1\right)+an_2\left(\left(1+\frac{1}{n-1}\right)^2-1\right)\\[2mm]
&\leq n+ c_3,
\end{align*}
where $c_3>0$ is a constant. Therefore, for sufficiently large $n$, there exists a constant $c_2>0$ such that
\begin{align}
\rho(H)-\rho(K_{q-1}\vee T_{p}(n-q+1))&\geq \frac{2f(k_q-1,k_q-1)}{n+ c_3}\nonumber\\[2mm]
&\geq \frac{2f(k_q-1,k_q-1)}{n}\frac{n}{n+c_3}\nonumber\\[2mm]
&\geq \frac{2f(k_q-1,k_q-1)}{n}\left(1-\frac{c_2}{n}\right).\label{Tnr-1}
\end{align}

\noindent{\bfseries Claim 3.} Let $H_{in}=\cup_{i=1}^{p}H[V_i\setminus W]$, and $H_{out}=(V(H)\setminus W,E(F)\setminus E(H\setminus W))$. Then
$e(H_{in})-e(H_{out})\leq f(k_q-1,k_q-1).$

We first prove that $H\setminus W$ is $S^{p+1}_{k_q}$-free. Suppose to the contrary that $H\setminus W$ contains a copy of  $S^{p+1}_{k_q}$, say $H_1$. Since $|W|=q-1$ and each vertex of $W$ has degree $n-1$, $H$ has a copy of    $S^{p+1}_{k_1,\cdots,k_{q-1}}$, say $H_2$, such that $V(H_1)\cap V(H_2)=\emptyset$. Thus $H_1\cup H_2$ is a copy of   $S^{p+1}_{k_1,\cdots,k_{q}}$ in $H$, and this is a contradiction.
From the definitions of  $H_{in}$ and $H_{out}$, we have  $e(H_{in})=\sum_{i=1}^{p}|E(H[V_i\setminus W])|$ and $e(H_{out})=\sum_{1\leq i<j\leq p}|V_i\setminus W||V_j\setminus W|-\sum_{1\leq i<j\leq p}|E(V_i\setminus W,V_j\setminus W)|$. By Lemma \ref{chen}, we only need to prove that for any $i\in [p]$ and  $v\in V_i\setminus W$,
\begin{eqnarray}
& &\sum_{j\neq i}\nu(H[V_j\setminus W])\leq k_q-1\ \ \  \mbox{and} \ \ \ \Delta(H[V_i\setminus W])\leq k_q-1,\label{eqn1}\\
& &d_{H[V_i\setminus W]}(v)+\sum_{j\neq i}\nu(H[N(v)\cap (V_j\setminus W)])\leq k_q-1.\label{eqn2}
\end{eqnarray}
Obviously  (\ref{eqn2}) implies  (\ref{eqn1}), so it is sufficient to prove (\ref{eqn2}).
We  prove  (\ref{eqn2}) by contradiction. Without loss of generality, suppose that there exists a vertex $u\in V_1\setminus W$ such that
$$d_{H[V_1\setminus W]}(u)+\sum_{j=2}^{p}\nu(H[N(u)\cap (V_j\setminus W)])\geq k_q.$$
Let $\{w_1w_2,\cdots,w_{2\ell-1}w_{2\ell}\}$ be an $\ell$-matching of $\cup_{j=2}^{p}H[N(u)\cap (V_j\setminus W)]$ and $u_1,\cdots,u_{k_q-\ell}\in V_1\setminus W$ be the neighbors of $u$. By Lemma \ref{Bi}, there exist $v_1,\cdots,v_{k_q-\ell}\in C_2$ such that $\{u,u_1,\cdots,u_{k_q-\ell},v_1,$ $\cdots,v_{k_q-\ell},w_1, w_2,\cdots,w_{2\ell}\}$ induce an $S^3_{k_q}$ of $H$.  For each $u_iv_i \ (1\leq i\leq k_q-\ell)$, there exist $p-2$ vertices $t_3\in C_3,t_4\in C_4,\cdots,t_{p}\in C_{p}$ such that $u,u_i,v_i,t_3,t_4,\cdots,t_{p}$ induce a $K_{p+1}$ of $H$. For any $ w_{i-1}w_i\in \{w_1w_2,\cdots,w_{2\ell-1}w_{2\ell}\}$, without loss of generality, suppose that $w_{i-1}w_i\in E(H[V_2\setminus W])$, then there exist $p-2$ vertices $z_3\in C_3,z_4\in C_4,\cdots,z_{p}\in C_{p}$ such that $u,w_{i-1},w_i,z_3,z_4,\cdots,z_{p}$ induce a $K_{p+1}$ of $H$. Thus $H\setminus W$ contains a copy of $S^{p+1}_{k_q}$, which contradicts that $H\setminus W$ is $S^{p+1}_{k_q}$-free.

According to the definitions of $H_{in}$, $H_{out}$ and $F$, we have $e(H)=e(K_{q-1}\vee F)+e(H_{in})-e(H_{out})$, and $x_v\leq x_{v_0}$ for any $v\in V(H_{in})$.
By Lemma \ref{Bi}, each vertex of $C_i$ is adjacent to all vertices of $V(H)\setminus V_i$ for any $i\in [p]$. Thus $$e(H_{out})\leq \sum_{1\leq i<j\leq p}|B_i||B_j|< \binom{p}{2}(2k_q^2)^2\leq 2k_q^4p^2.$$
By Lemma \ref{eigenvector}, we have $$\mathbf{x}^{\mathrm{T}}\mathbf{x}\geq n(x_{v_0}-\frac{100k_1pq}{n})^2\geq nx_{v_0}^2-200k_1pqx_{v_0}.$$  Therefore,
\begin{align}
\rho(H)-\rho(K_{q-1}\vee F)&\leq  \frac{\mathbf{x}^\mathrm{T}(A(H)-A(K_{q-1}\vee F))\mathbf{x}}{\mathbf{x}^{\mathrm{T}}\mathbf{x}}\nonumber\\[2mm]
&= \frac{2\sum_{ij\in E(H_{in})}\mathbf{x}_i\mathbf{x}_j}{\mathbf{x}^{\mathrm{T}}\mathbf{x}}- \frac{2\sum_{ij\in E(H_{out})}\mathbf{x}_i\mathbf{x}_j}{\mathbf{x}^{\mathrm{T}}\mathbf{x}} \nonumber\\[2mm]
&\leq  \frac{2e(H_{in})x_{v_0}^2}{\mathbf{x}^{\mathrm{T}}\mathbf{x}}- \frac{2e(H_{out})(x_{v_0}-\frac{100k_1pq}{n})^2}{\mathbf{x}^{\mathrm{T}}\mathbf{x}} \nonumber\\[2mm]
&\leq \frac{2(e(H_{in})-e(H_{out}))x_{v_0}^2}{\mathbf{x}^{\mathrm{T}}\mathbf{x}}+\frac{2e(H_{out})\frac{200k_1pq}{n}}{\mathbf{x}^{\mathrm{T}}\mathbf{x}}\nonumber\\[2mm]
&\leq \frac{2f(k_q-1,k_q-1)x_{v_0}^2}{nx_{v_0}^2-200k_1pqx_{v_0}}+\frac{\frac{800k_1k^{4}_qp^3q}{n}}{nx_{v_0}^2-200k_1pqx_{v_0}}\nonumber\\[2mm]
&\leq  \frac{2f(k_q-1,k_q-1)}{n-\frac{200k_1p^2q}{p-2}}+\frac{800k_1k^{4}_qp^3q}{n^2x_{v_0}^2-200k_1pqn}\label{GC}
\end{align}

Combining (\ref{Tnr-1}) and (\ref{GC}), there exists a constant $c_4>0$ such that
\begin{eqnarray*}
& &\rho(K_{q-1}\vee T_{p}(n-q+1))-\rho(K_{k-1}\vee F)\\[2mm]
&\leq & \frac{2f(k_q-1,k_q-1)}{n-\frac{200k_1p^2q}{p-2}}+\frac{800k_1k^{4}_qp^3q}{n^2x_{v_0}^2-200k_1pqn}-\frac{2f(k_q-1,k_q-1)}{n}\left(1-\frac{c_2}{n}\right)\\[2mm]
&< & \frac{2f(k_q-1,k_q-1)\frac{200k_1p^2q}{p-2}}{n(n-\frac{200k_1p^2q}{p-2})}+\frac{800k_1k^{4}_qp^3q}{n^2x_{v_0}^2-200k_1pqn}+\frac{2f(k_q-1,k_q-1)c_2}{n^2}\\[2mm]
&< & \frac{c_4}{n^2}.
\end{eqnarray*}
Combining with Claim 1, we have
 \begin{align*}
 \frac{c_1}{n}\leq \rho(K_{q-1}\vee T_{p}(n-q+1))-\rho(K_{q-1}\vee F) \leq \frac{c_4}{n^2},
 \end{align*}
 which is a contradiction when $n$ is sufficiently large. Thus $\left||V_i\setminus W|-|V_j\setminus W|\right|\leq 1$ for any $1\leq i<j\leq p$.

\qed

\medskip
\noindent{\bfseries Proof of Theorem \ref{main}.} Now we  prove that $e(H)=\mathrm{ex}(n,S^{p+1}_{k_1,\cdots,k_q})$. Otherwise, we  assume that $e(H)\leq \mathrm{ex}(n,S^{p+1}_{k_1,\cdots,k_q})-1$. Let $H'$ be an $S^{p+1}_{k_1,\cdots,k_q}$-free graph with $e(H')=\mathrm{ex}(n,S^{p+1}_{k_1,\cdots,k_q})$.
By Theorem \ref{extremal}, $H'$ can be obtained from $K_{q-1}\vee T_{p}(n-q+1)$ by embedding some graph in $T_{p}(n-q+1)$. By Lemma \ref{balance}, we may assume that $V_1\cup\cdots \cup V_{p}$ is a vertex partion of $H'$. Let $E_1=E(H)\setminus E(H')$, $E_2=E(H')\setminus E(H)$,  then $E(H')=(E(H)\cup E_2)\setminus E_1$, and
\[
|E(H)\cap E(H')|+|E_1|=e(H)<e(H')=|E(H)\cap E(H')|+|E_2|,
\]
which implies that
$|E_2|\geq |E_1|+1$. Furthermore, by Lemma \ref{Bi}, we have
\begin{eqnarray}
|E_2|=e(H')-e(H)\leq f(k_q-1,k_q-1)+\sum\limits_{1\leq i<j\leq p}|B_i||B_j|< k^2_q +\binom{p}{2}(2k_q^2)^2< 3k_q^4p^2.\label{eqn3}
\end{eqnarray}
Therefore,
\begin{align*}
\rho(H')-\rho(H)&\geq \frac{\mathbf{x}^{\mathrm{T}}(A(H')-A(H))\mathbf{x}}{\mathbf{x}^{\mathrm{T}}\mathbf{x}}\\[2mm]
&= \frac{2\sum_{ij\in E_2}\mathbf{x}_i\mathbf{x}_j}{\mathbf{x}^{\mathrm{T}}\mathbf{x}}-\frac{2\sum_{ij\in E_1}\mathbf{x}_i\mathbf{x}_j}{\mathbf{x}^{\mathrm{T}}\mathbf{x}}\\[2mm]
&=\frac{2}{\mathbf{x}^{\mathrm{T}}\mathbf{x}}\Big(\sum_{ij\in E_2}\mathbf{x}_i\mathbf{x}_j-\sum_{ij\in E_1}\mathbf{x}_i\mathbf{x}_j\Big)\\[2mm]
&\geq \frac{2}{\mathbf{x}^{\mathrm{T}}\mathbf{x}}\Big(|E_2|(x_{v_0}-\frac{100k_1pq}{n})^2- |E_1|x^2_{v_0}\Big)\\[2mm]
&\geq  \frac{2}{\mathbf{x}^{\mathrm{T}}\mathbf{x}}\Big(|E_2|x^2_{v_0}-\frac{200k_1pqx_{v_0}}{n}|E_2|- |E_1|x^2_{v_0}\Big)\\[2mm]
&\geq \frac{2x_{v_0}}{\mathbf{x}^{\mathrm{T}}\mathbf{x}}\Big(x^2_{v_0}-\frac{200k_1pq}{n}|E_2|\Big)\\[2mm]
&\geq \frac{2x_{v_0}}{\mathbf{x}^{\mathrm{T}}\mathbf{x}}\Big((1-\frac{2}{p})^2-\frac{200k_1pq}{n} 3k_q^4p^2\Big)\\[2mm]
&> 0,
\end{align*}
which contradicts the assumption that $H$ has the maximum spectral radius among all $n$-vertex $S^{p+1}_{k_1,\cdots,k_q}$-free graphs. Hence $e(H)=\mathrm{ex}(n,S^{p+1}_{k_1,\cdots,k_q})$.
\qed

\end{document}